\documentclass{article}

\usepackage{amsmath}
\usepackage{amsmath,amssymb, amsthm, latexsym, mathrsfs,amsfonts, epsfig, color,graphicx,}

\renewcommand{\a}{\alpha}

\newtheorem{theorem}{Theorem}

\newtheorem{remark}{Remark}

\newtheorem{definition}{Definition}

\begin{document}
\title{Existence and Uniqueness of a Fractional Fokker-Planck Equation}
\author{Li Lin \footnote{E-mail address: linli@hust.edu.cn} \\
\small{ Center for Mathematical Sciences, Huazhong University of Science and Technology} \\
\small {Wuhan, 430074, China}\\}
\maketitle
\begin{abstract}
Stochastic differential equations with L\'{e}vy motion arise the mathematical models for various phenomenon in geophysical and biochemical sciences.
The Fokker-Planck equation for such a stochastic differential equations is a nonlocal partial differential equations. We prove the existence and uniqueness of the weak solution for this equation.
\end{abstract}
keywords: Fractional Laplacian operator, Fokker-Planck equation, Lax- Milgram theorem, Existence and  uniqueness

\section{introduction}
For a system described by a stochastic differential equation with a $\alpha-$stable L\'{e}vy motion $L_t^{\a}$ for $\a \in(0,2),$
$$dX_t=b(X_t)dt+adB_t+dL_t^{\a},$$
the corresponding Fokker-Planck equation \cite{D} contains a nonlocal Laplacian operator,
$$u_t=a\Delta u-(-\Delta)^{\a/2}u-div(b(x)\cdot u),$$
where $b$ is a two-dimensional vector function and  $-(-\Delta)^{\a/2}$ is a nonlocal Laplacian operator defined by
$$(-\Delta)^{\a/2}f(x)=\int_{\mathbb{R}^d/\{0\}}\frac{f(x)-f(y)}{|x-y|^{d+\a}}dy.$$

Fractional partial differential equations  arise the mathematical models for various phenomenon in physics and biology.
Such as anomalous diffusion of particles \cite{MM} and  the cell density evolution in certain biological processes \cite{JD}. Volume constraints are natural extensions to the fractional case of boundary conditions for differential equations \cite{QD}.

The Fokker-Planck equations have been used in the modeling of many physical phenomena, in particular, for the description of the evolution of plasmas \cite{HR}. Moreover, there are some settings in which particles may have long jumps\cite{DS}. Questions such as existence of solutions, hydrodynamic limits, and long-time behavior for the Fokker-Planck system has been extensively studied by many authors\cite{FB,NE,TG} .

There are few works dealing with the existence and uniqueness of the Fokker-Planck equations. Wei and Tian \cite{Wei}
obtain the existence and uniqueness of weak $L^p$ solution on the whole space. AcevesSanchez and Cesbron \cite{PA} studied the nonlocal Fokker-Planck problem on $\mathbb{R}^d$ in fractional Sobolev space with the drift term is a direct proportion function.  In \cite{He}, the weak solution is just considered in the Sobolev space $H^1.$ Moreover, there is no proof of the well-posedness for equation  with Laplace operator.

This paper is devoted to the study of the nonlocal Fokker-Planck equation on a bounded interval in $\mathbb{R}^2.$ There are two main results: For $a\geq0,$ we obtain the existence and uniqueness in the Sobolev space $H^1$ with respect to spatial variables; for $a=0,$ we reach the conclusion in a fractional sobolev space $H^{\a/2}$ with respect to spatial variables. We can see that, This is more accurate.
\section{Well-posedness of  Fokker-Planck equation}
In this section, we give some function spaces and obtain  the existence and uniqueness of
the weak solution of a fractional Fokker-Planck equation on a bounded domain in $\mathbb{R}^2.$
\subsection{Function spaces}
Denote $Q_T=D\times(0,T).$ Let $\partial_lQ_T$ and $\partial_pQ_T$ be the lateral boundary $D^c\times(0,T)$ and the parabolic boundary $\partial_lQ_T\cup \{(x,t); x\in \bar{D},t=0\}$ of $Q_T.$
Denote by $\mathring{W}, \dot{W}$ the set of all functions in $W$ vanish on $\partial_pQ_T$ and $\partial_lQ_T$ respectively.

For the integer $\beta, r \in \{0,1\},$ the set
$$\{u; \nabla^\beta u, \partial^r_t u \in L^2(Q_T)\}$$
endowed with the norm
$$||u||^2_{W_2^{1,1}(Q_T)}=\int\int_{Q_T}\sum\limits_{\beta\in\{0,1\}}|\nabla^\beta u|^2+\sum\limits_{r\in\{0,1\}}|\partial^r_t u|^2dxdt$$
is denoted by $W_2^{1,1}(Q_T).$
Denote
$$V(Q_T):=\{u\in \dot{W}_2^{1,1}(Q_T), \nabla u_t \in L^2(Q_T;\mathbb{R})\},$$
and define the inner product as
$$(u,v)_{V(Q_T)}=(u,v)_{W_2^{1,1}(Q_T)}+(\nabla u_t,\nabla v_t)_{L^2(Q_T)}.$$

\subsection{Main result}
We consider the following nonlocal Fokker-Planck equation
\begin{equation}{\label{1}}
  \begin{cases}
    u_t=a\Delta u-(-\Delta)^{\a/2}u-div(b(x)\cdot u)\qquad x\in D\subset \mathbb{R}^2,\;t\in(0,T),\\
    u|_{D^c}=0,\\
    u(x,0)=u_0(x),
  \end{cases}
\end{equation}
where $D=(0,1)$ is an bounded domain in $\mathbb{R}^2,$ and $D^c$ is the complement of D. We will prove the existence and uniqueness of the solution to the equation (\ref{1}).

\begin{remark}
  We define nonlocal divergence operator $\mathcal{D}$ on $\beta$ as
$$\mathcal{D}(\beta)(x):=\int_{\mathbb{R}^d}(\beta(x,z)+\beta(z,x))\cdot\gamma(x,z)dz  \qquad \text{for}\; x\in D.$$
For a function  $\phi(x)$, the adjoint operator $\mathcal{D^{*}}$ corresponding to $\mathcal{D}$ is the operator whose action on $\phi$ is given by
$$\mathcal{D^{*}}(\phi)(x,z)=-(\phi(z)-\phi(x))\gamma(x,z) \qquad \text{for} \; x,z\in D.$$
Here we take $\gamma(x,z)=\frac{1}{\sqrt{2}}(z-x)\frac{1}{|z-x|^{\frac{d+2+\a}{2}}}.$
we have
$$\mathcal{D}\mathcal{D^{*}}=-(-\Delta)^{\a/2}.$$
\end{remark}

\begin{definition}
  Consider $u_0$ in $L^2(D).$ We say that $u$ is a weak solution of equation (\ref{1}), if for any $\varphi\in C_c^\infty([0,T)\times D)$,
\begin{equation}
  \int\int_{Q_T}u(\partial_t \varphi +a\Delta \varphi-(-\Delta)^{\a/2}\varphi+b(x)\cdot\nabla\varphi)dxdt+ \int\int_{Q_T}u_0(x)\varphi(0,x)dx=0
\end{equation}
 \end{definition}

\begin{theorem}\label{t1}
Consider $u_0\in L^2(D)$,
\begin{enumerate}
  \item If $a=0$ and $div b(x)\geq0$, there exists a unique weak solution and this solution satisfies
$$f\in \mathcal{X}:=\{f: f\in L^2(Q_T), \frac{|f(t,x)-f(t,y)|}{|x-y|^{\frac{2+\a}{2}}}\in L^2(Q_T\times\mathbb{R}^2)\}.$$
  \item For $a\geq0$, there exist a unique weak solution in $\mathring{W}_2^{1,1}(Q_T).$
\end{enumerate}
\begin{remark}
  Note that this definition of $\mathcal{X}$ is equivalent to saying that it is the set of functions which are in $L^2([0,T))$ with respect to time  and in $H^{\alpha/2}(D)$ with respect to space.
\end{remark}
\begin{proof}
  $1.$ We consider the Hilbert space $\mathcal{X}$ provided with the norm
  $$||u||_{\mathcal{X}}^2=||u||^2_{L^2(Q_T)}+||D^*u||^2_{L^2(Q_T\times\mathbb{R})}.$$
  Let us denote $\mathcal{T}$ the operator, given by
  $$\mathcal{T}u=\partial_tu+div(b(x)\cdot u).$$
  Moreover, we define the Hilbert space $\mathcal{Y}$ as
  $$\mathcal{Y}=\{f\in \mathcal{X}: \mathcal{T}u\in \mathcal{X}^{'}\},$$
  where $\mathcal{X}^{'}$ is the dual of $\mathcal{X}.$
  From the fact that the space $C_c^\infty(Q_T)$ is a subspace of $\mathcal{X}$ with a  continuous injection, we define the pre-Hilbertian norm:
  $$|\varphi|_{C_c^\infty(Q_T)}=||\varphi||_{\mathcal{X}}^2+\frac{1}{2}||\varphi(0,x)||^2_{L^2(D)}.$$
  Now, we introduce the bilinear form $a: \mathcal{X}\times C_c^\infty(Q_T)\rightarrow \mathbb{R}$ as
  \begin{equation}\label{2}
  a(u,\varphi)=\int\int_{Q_T}-u\varphi_t+(\mathcal{D}^*u,\mathcal{D}^*\varphi)_{L^2(\mathbb{R})}-b(x)u\cdot\nabla \varphi dxdt.
  \end{equation}
  and $$L(\varphi)=-\int_{D}u_0(x)\varphi(0,x)dx.$$
  Then find a solution $u$ in $\mathcal{X}$ of $(\ref{1})$ is equivalent to finding a solution $u$ in $\mathcal{X}$ of $a(u,\varphi)=L(\varphi)$ for any $\varphi\in C_c^\infty(Q_T)$. First $a(u,\varphi)$ is continuous. Next we will obtain the coercivity of $a$.
  \begin{equation}\label{5}
  \begin{split}
  a(\varphi,\varphi)&=\int\int_{Q_T}-\varphi_t\varphi+(\mathcal{D}^*\varphi,\mathcal{D}^*\varphi)_{L^2{(\mathbb{R})}}-b(x)\varphi\cdot\nabla\varphi dxdt\\
  &=\frac{1}{2}\int_{D}\varphi^2(0,x)dx+\int\int_{Q_T}(\mathcal{D}^*\varphi,\mathcal{D}^*\varphi)_{L^2{(\mathbb{R})}}dxdt+\frac{1}{2}\int_{D}div b(x)\varphi^2dxdt.
  \end{split}
  \end{equation}

  Then, there exists a positive constant $\delta$ such that  $a(\varphi,\varphi)\geq \delta|\varphi|_{C_c^\infty(Q_T)}.$
  Thus the Lax-Milgram theorem implies the existence of $u$ in $\mathcal{X}$ satisfying equation (\ref{2}). That yields existence of a solution of equation (\ref{1}) with $a=0.$ For $u\in \mathcal{X},$ the linear bounded operator $\mathcal{T}$  maps $u\in \mathcal{T}$ to $-(-\Delta)^{\a/2}u \in \mathcal{X}^{'}$, hence
  the weak solution $u$ is in $\mathcal{Y}.$

  Since the equation (\ref{1}) is linear, to show the uniqueness, it is  enough to show the unique solution with zero initial is the function $u\equiv 0.$ Let $u$ be a solution of this problem on $\mathcal{Y}.$
  Through integration by parts we have
  \begin{equation*}
  \begin{split}
  2(\mathcal{T}u,u)_{\mathcal{X},\mathcal{X}^{'}}&=(\partial_tu+div(b(x)\cdot u),u)_{\mathcal{X},\mathcal{X}^{'}}\\
  &=\int_{D}u^2(T,x)dx+\int_{Q_T}divb(x)u^2dxdt\geq0.
  \end{split}
  \end{equation*}
  On the other hand, since $u$ satisfies equation (\ref{1}), $\mathcal{T}u=-(-\Delta)^{\a/2}$ in the weak sense, then
  $$(\mathcal{T}u,u)_{\mathcal{X},\mathcal{X}^{'}}=-(\mathcal{D}^*u,\mathcal{D}^*u)\leq0.$$
  That means $u\equiv0 $ a.e. on $Q_T$. The solution is unique.

$2.$ For the first part we set $a>0.$
To show the unique of the equation (\ref{1}), it is  enough to show the unique solution with zero initial is the function $w\equiv 0.$
\begin{equation}
  \begin{cases}
    w_t=a\Delta w-(-\Delta)^{\a/2}w-div(b(x)\cdot w) \qquad x\in D\subset \mathbb{R}^2,\\
    w|_{D^c}=0,\\
    w(x,0)=0.
  \end{cases}
\end{equation}

Set $B[u,v,t]=\int_Da\nabla u\cdot\nabla v+(D^*u,D^*v)_{L^2(\mathbb{R})}-b(x)\cdot u\cdot\nabla udx,$
then $$\frac{d}{dt}(\frac{1}{2}||w||^2_{L^2(D)})+B[w,w,t]=(w,w')+B[w,w,t]=0.$$
This means
\begin{equation}
\begin{split}
  a\int_{D}|\nabla w|^2dx&+(D^*w,D^*w)_{L^2(D\times\mathbb{R})}\\
  &=B[w,w,t]+\int_D b(x)\cdot w\cdot\nabla wdx\\
  &\leq B[w,w,t]+\varepsilon\int_{D}|\nabla w|^2dx+\frac{C}{4\varepsilon}\int_{D}|w|^2dx.
\end{split}
\end{equation}
From Poincar\'{e} inequality $$||w||_{L^2(D)}\leq C_1||D^*w||_{L^2(D\times\mathbb{R})}.$$
Then choose $\varepsilon$ small enough, then
$$\sigma||w||_{H^1}+\beta||w||_{H^{\a/2}}\leq B[w,w,t]+\gamma ||w||^2_{L^2(D)},$$
for positive constant $\sigma, \beta,$ and $\gamma\geq 0.$
We can see that $$B[w,w,t]\geq -\gamma ||w||^2_{L^2(D)}.$$
By Gronwall inequality, we obtain
$$w=0.$$ The solution is unique.

Next, we will prove the existence of the problem. Without loss of generality, we just consider the case $u_0=0.$
$$a(u,v)=\int\int_{Q_T}(u_tv_t+a\nabla u\nabla v_t+( D^*u, D^*v_t)_{L^2(\mathbb{R})}+ div(b(x)\cdot u)v_t)e^{-\theta t}dxdt, \theta>0.$$
Then $$|a(u,v)|\leq ||u||_{\dot{W}^{1,1}(Q_T)}||v||_{V(Q_T)}.$$
set $$a(v,v):=A+B+C+D$$
For $v\in V(Q_T)$,
\begin{equation}
  \begin{split}
   B&=a\int\int_{Q_T}\nabla v\cdot\nabla v_te^{-\theta t}dxdt
 =a\int\int_{Q_T}\frac{1}{2}\frac{\partial}{\partial_t}|\nabla v|^2e^{-\theta t}dxdt\\
 &=a\int\int_{Q_T}\frac{1}{2}\frac{\partial}{\partial_t}(|\nabla v|^2e^{-\theta t})dxdt+\frac{\theta}{2}\int\int_{Q_T}e^{-\theta t}|\nabla v|^2dxdt\\
 &=\frac{a}{2}e^{-\theta T}\int_{D}|\nabla v|^2|_{t=T}dx-\frac{a}{2}\int_{D}\gamma|\nabla v|^2|_{t=0}dx+\frac{a\theta}{2}\int\int_{Q_T}e^{-\theta t}|\nabla v|^2dxdt\\
  \end{split}
\end{equation}
Since $v\in V(Q_T)$, then $\gamma|\nabla v|^2_{t=0}=0,$ combine with Poincar\'{e} inequality
\begin{equation}
  \begin{split}
  B&\geq \frac{\theta}{2}\int\int_{Q_T}e^{-\theta t}|\nabla v|^2dxdt\\
  &\geq\frac{a\theta}{4}\int\int_{Q_T}e^{-\theta t}|\nabla v|^2dxdt+\frac{a\theta}{4\mu}\int\int_{Q_T}e^{-\theta t}v^2dxdt.
 \end{split}
\end{equation}

Similarly,
$$C\geq\frac{\theta}{2}\int\int_{Q_T}||\mathcal{D}^*v||_{L^2(\mathbb{R})}^2dxdt.$$
On the other hand, we have
\begin{equation}
  \begin{split}
  D&=\int\int_{Q_T} div(b(x)\cdot v)v_te^{-\theta t}dxdt\\
  &=\int\int_{Q_T}divb(x)\cdot v\cdot v_te^{-\theta t}dxdt+\int\int_{Q_T}b(x)\cdot\nabla v\cdot v_te^{-\theta t}dxdt\\
  &:=E+F,
   \end{split}
\end{equation}
where
$$|E|\leq \varepsilon\int\int_{Q_T}v_t^2e^{-\theta t}dxdt+\frac{C_2}{\varepsilon}\int_{Q_T}v^2e^{-\theta t}dxdt.$$
$$|F|\leq \varepsilon\int\int_{Q_T}v_t^2e^{-\theta t}dxdt+\frac{C_3}{\varepsilon}\int\int_{Q_T}|\nabla v|^2e^{-\theta t}dxdt.$$
Then we have
\begin{equation}
  \begin{split}
  a(v,v)&\geq(1-2\varepsilon)\int_{Q_T}v_t^2e^{-\theta t}dxdt+\frac{\theta}{2}\int_{Q_T}||D^*v||^2_{L^2(\mathbb{R})}e^{-\theta t}dxdt\\
  &(\frac{a\theta}{4}-\frac{C_4}{\varepsilon})\int_{Q_T}|\nabla v|^2e^{-\theta t}dxdt+(\frac{a\theta}{4\mu}-\frac{C_5}{\varepsilon})\int_{Q_T}v^2e^{-\theta t}dxdt
   \end{split}
\end{equation}
Choose $\varepsilon$ small enough and $\theta>0$ large enough, we have
$$ a(v,v)\geq\delta||v||_{W^{1,1}(Q_T)},$$
where $\delta$ is a positive constant.
Then the Lax-Milgram theorem implies the existence of equation (\ref{1}).

For $a=0,$ it is similar to the work that He and Duan did in \cite{He}.
\end{proof}
\end{theorem}

\section{Example}
Langevin equation provide models of a diffusing particle. We consider the following system of stochastic differential equations
\begin{equation}\label{3}
  \begin{cases}
    dx_t=v_tdt,\\
    mdv_t=-\gamma v_tdt+\sigma dL_t^{\a},
  \end{cases}
\end{equation}
where $m$ is the mass of the particle, $\gamma$ and $\sigma$ are the dissipation and
diffusion coefficient, respectively, and $L_t^{\a}$ is a L\'{e}vy process with generator
$-(-\Delta)^{\a/2}.$
Then the corresponding L\'{e}vy Fokker-Planck equation is
\begin{equation}\label{6}
  \begin{cases}
 u_t+v\cdot \nabla_xu=-\frac{\sigma}{m}(-\Delta_v)^{\a/2}u+\frac{\gamma}{m}div_v(v\cdot u)\qquad x,v\in D\subset \mathbb{R}^2,\;t\in(0,T),\\
    u|_{D^c}=0,\\
    u(x,0)=u_0(x)
  \end{cases}
\end{equation}
We can define the bilinear form as (\ref{2})
$$a_1(u,\varphi)=\int\int\int_{Q_T\times D}-u\varphi_t-vu\cdot\nabla_x\varphi+\frac{\sigma}{m}(\mathcal{D}_v^*u,\mathcal{D}_v^*\varphi)_{L^2(\mathbb{R})}+\frac{\gamma}{m}vu\cdot\nabla_v \varphi dxdtdv.$$
and from the fact $div_v(v)>0,$ and
$$\int\int\int_{Q_T\times D} v \varphi \cdot\nabla_x \varphi dxdtdv=0,$$
  Then, there exists a positive constant $\delta_1$ such that  $a(\varphi,\varphi)\geq \delta_1|\varphi|_{C_c^\infty(Q_T\times D)}.$
we verified that the equation (\ref{6}) has a unique weak solution in $\mathcal{X}_1$ by Theorem \ref{t1},
where, $$\mathcal{X}_1:=\{f: f\in L^2(Q_T\times D), \frac{|f(t,x,v)-f(t,x,w)|}{|v-w|^{\frac{2+\a}{2}}}\in L^2(Q_T\times D\times \mathbb{R}^2)\}.$$
Clearly, the result of the example agree with the theoretical finding in this study.

\end{document}